\documentclass[11pt,a4paper, twoside]{article}
\usepackage[latin9]{inputenc}
\usepackage{amsmath,amssymb,amsthm,amscd,hyperref}
\usepackage[]{fontenc}
\usepackage{xy}
\usepackage{enumerate}
\usepackage{graphicx}
\usepackage[capitalise]{cleveref}

\newif\ifpdf\ifx\pdfoutput\undefined\pdffalse\else\pdfoutput=1\pdftrue\fi

\xyoption{all} \numberwithin{equation}{section}
\newtheorem{theorem}[equation]{Theorem}

\newtheorem{proposition}[equation]{Proposition}
\newtheorem{lemma}[equation]{Lemma}

\newtheorem{conjecture}[equation]{Conjecture}
\theoremstyle{definition}

\newtheorem{definition}[equation]{Definition}
\newtheorem{notation}[equation]{Notation}

\newtheorem{remark}[equation]{Remark}

\newcommand{\nc}{\newcommand}
\nc{\cH}{\mathcal{H}} \nc{\cG}{\mathcal{G}}
\nc{\cC}{\mathcal{C}}
\nc{\cO}{\mathcal{O}}
\nc{\cI}{\mathcal{I}}
\nc{\cA}{\mathcal{A}}
\nc{\cB}{\mathcal{B}} \nc{\cY}{\mathcal{Y}} \nc{\cK}{\mathcal{K}}
\nc{\cX}{\mathcal{X}} \nc{\cS}{\mathcal{S}} \nc{\cE}{\mathcal{E}}
\nc{\cF}{\mathcal{F}} \nc{\cZ}{\mathcal{Z}} \nc{\cQ}{\mathcal{Q}}
\nc{\cN}{\mathcal{N}} \nc{\cP}{\mathcal{P}} \nc{\cL}{\mathcal{L}}
\nc{\cM}{\mathcal{M}} \nc{\cR}{\mathcal{R}} \nc{\cT}{\mathcal{T}}
\nc{\cW}{\mathcal{W}} \nc{\cU}{\mathcal{U}} \nc{\cD}{\mathcal{D}}
\nc{\cJ}{\mathcal{J}} \nc{\cV}{\mathcal{V}}
\nc{\fr}{{\rightarrow}}
\nc{\rd}{red.deg}
\newcommand{\Q}{\mathbb{Q}}
\newcommand{\Z}{\mathbb{Z}}

\newcommand{\pr}{\mathbb P}

\newcommand{\sym}{\mbox{\upshape{Sym}}}

\newcommand{\rk}{\mbox{\upshape{rank}}}

\newcommand{\cliff}{\mbox{\upshape{Cliff}}}
\newcommand{\alb}{\mbox{\upshape{alb}}}
\newcommand{\Alb}{\mbox{\upshape{Alb}}}
\newcommand{\pic}{\mbox{\upshape{Pic}}}

\title{Fibred surfaces and their unitary rank}
\author{Lidia Stoppino}
\date{\small{Dedicated to the Memory of Gang Xiao}}

\begin{document}
\maketitle
\begin{abstract}
Let $f\colon S\to B$ a complex fibred surface with fibres of genus $g\geq 2$. Let $u_f$ be its unitary rank, i.e., the rank of the maximal unitary summand of the Hodge bundle $f_*\omega_f$. 
We prove many new slope inequalities involving $u_f$ and some other invariants of the fibration. As applications: 
\begin{itemize}
\item[(1)]we prove a new Xiao-type bound on $u_f$ with respect to $g$ for non-isotrivial fibrations:
\[
u_f< g\frac{5g-2}{6g-3}.
\] 
In particular this implies that if $f$ is not locally trivial and $u_f=g-1$ is maximal, then $g\leq 6$; 
\item[(2)] we prove a result in the direction of the Coleman-Oort conjecture: a new constraint on the rank of the $(-1,0)$ part of the maximal unitary Higgs subbundle of a curve generically contained in the Torelli locus.
\end{itemize}
\end{abstract}


\section{Introduction}\label{intro}

This paper deals with the geography of fibred surfaces and some applications. We work over the complex field $\mathbb C$.
The main results can be divided in four sections, which we list in order of appearance.
\subsection{New slope inequalities}
Let $f\colon S\to B$ be a relatively minimal fibred surface. 
The relative canonical line bundle is $\omega_f:=\omega_S\otimes f^*\omega_B^{-1}$; we call $K_f$ any associated divisor. 
Let $F$ be a general fibre of $f$.  Let $\chi_f:=\chi(\cO_S)-\chi(\cO_F)\chi(\cO_B)$ be the relative  Euler characteristic (see \Cref{sec: preliminaries}). 
 
We prove new {\em slope inequalities,} i.e., bounds of the form 
\[K_f^2\geq \alpha\chi_f,\]
where $\alpha$ is a function of some invariants of the fibration. In particular, in our case $\alpha$ will be a function depending on the genus $g$ of the general fibre, on the unitary rank $u_f$, and in some cases on the Clifford index of the general fibre $c_f$.
The most significant inequalities we prove are the following:

\begin{theorem}[\Cref{teo: konnoplus}, \Cref{teo: eurekina}]\label{teo: principale}
Let $f\colon S\to B$ be a relatively minimal locally non-trivial fibred surface of genus $g\geq 2$.
\begin{itemize}
\item The following inequality holds:
\begin{equation} \label{eq: konnoplus1}
K_f^2 \geq \frac{4g(g-1)}{(2g-1)(g-u_f)}\chi_f.
\end{equation}
Moreover, the equality holds in \eqref{eq: konnoplus1} if and only if $u_f=g-1$, the general fibre $F$ is trigonal and there exists a section $\Gamma$ of $f$ such that $K_f $ is numerically equivalent to  $(2g-2)\Gamma+ \chi_f F$.
\item Assume  that the ample summand  $\cA$  in the second Fujita decomposition \eqref{second fujita} of $f_*\omega_f$ is semistable. 
Then the following inequality holds:
\begin{equation}\label{eq: eureka1}
K_f^2\geq \left( \frac{4(g-1)-2u_f}{(g-u_f)}+ \frac{2u_f}{(2u_f+1)(g-u_f)}\right)\chi_f.
\end{equation}
Moreover, if the equality holds in \eqref{eq: eureka1}, then necessarily
$u_f=g-1$ and $F$ is trigonal.
\item  Assume  that $\cA$ is semistable and that $u_f\leq c_f$. Then the following inequality holds:
\begin{equation}\label{eq: eurekabis1}
K_f^2\geq \left( \frac{4(g-1)-u_f}{(g-u_f)}+ \frac{u_f}{(u_f+1)(g-u_f)}\right)\chi_f.
\end{equation}

\end{itemize}
\end{theorem}
\begin{remark}
The first inequality \eqref{eq: konnoplus1} is an extension of a result of Konno in \cite{konno}, where he proved the 
same inequality with the relative irregularity $q_f$ instead of $u_f$. 
The second and third inequalities are a sharpening of an inequality due to the author with Barja in \cite{SB} and of inequalities obtained by the author with Riva in \cite{R-S}.
\end{remark}

\subsection{A new Xiao type inequality for the unitary rank}
We extend the bound obtained by Konno in \cite{konno}  to the unitary rank. Konno's bound was itself a sharpening of the famous  Xiao's  bound \cite{xiao-irregular}. 
Moreover, our bound is strict.
In particular, we obtain that:
\begin{theorem}[\Cref{thm: bound}, \Cref{prop: estremale}]
For a relatively minimal  non-isotrivial fibred surface of genus $g$ and  unitary rank $u_f$, it always holds 
\begin{equation}\label{eq: buh}
u_f< g\frac{5g-2}{6g-3}.
\end{equation}
Moreover, the following hold:
\begin{enumerate} 
\item if the unitary rank is maximal $u_f=g-1$, then necessarily $g\leq 6$;
\item if the unitary rank is $u_f=g-2$, then necessarily $g\leq 11$.
\end{enumerate}
\end{theorem}
\noindent The bound \eqref{eq: buh} is, to our knowledge, the sharpest known bound for $u_f$; see \Cref{sec: preliminaries} for the known results and conjectures.
\subsection{Results on the case with high unitary rank} 
Pirola proved in \cite{pietroBUMI} that, for a fibration such that $u_f=g-1$, up to an \'etale base change the relative irregularity coincides with the unitary rank $q_f=u_f$.
Also the following result is essentially due to Pirola  (in a private communication with the author) and to the already known results \cite{BGN}, \cite{GTS}.
\begin{theorem}[\Cref{prop: estremale2}, \Cref{teo: estremale}]
Let $f\colon S\to B$ a relatively minimal locally non-trivial fibred surface of genus $g\geq 4$ such that $u_f= g-1$. Then $f$ is a non-isotrivial trigonal fibration and 
its smooth fibres $F_t$ are a covering of a non-isotrivial family of elliptic curves $E_t$. Moreover, $f$ is not Kodaira, i.e., $K_f^2<12\chi_f$.
\end{theorem}
\subsection{Applications to the Coleman-Oort conjecture}
We apply the inequalities of \Cref{teo: principale} to the study of Shimura varieties on the moduli space $\cA_g$ of principally polarized abelian varieties of genus $g$.
We work in the framework established by Viehweg and Zuo in \cite{VZ} and developed by  Chen, Lu and Zuo in \cite{CLZ} and in \cite{LZ4}. 
We obtain a result that goes in the direction of the Coleman-Oort conjecture. See  \Cref{sec: CO} for the notation and some explanation.
\begin{theorem}[\Cref{teo: CO}]
Let $C \subset \cA_g$ be a curve with Higgs bundle decomposition $\cE_C = \cA_C \oplus \cU_C$, where $\cU_C$ is the maximal unitary Higgs subbundle.

\smallskip
If
$\rk \,\cU_C^{-1,0 }\geq  g \frac{5g-2}{6g-3}$, then $C$ is not contained generically in the Torelli locus.

Moreover, if $C$ is a Shimura curve such that $\rk \,\cU_C^{-1,0} \geq  \frac{4g-1+\sqrt{16g^2-36g+21}}{10}$ then $C$ is not contained generically in the Torelli locus.
\end{theorem}
This is a sharpening of a theorem of Chen, Lu and Zuo \cite[Thm 1.1.2]{CLZ}: see  \Cref{rem: confronto LZ}.

\bigskip

A few words about the proofs. We work in the framework of the celebrated Xiao's method, with at least two  new key ingredients:
first we improve a result due to Konno in \cite{konno} that gives  an original construction of a nef divisor on the surface (see also \Cref{rem: novita'}).
Here we push Konno's idea to his full capacity, connecting the results with the unitary rank of the fibration, while in his paper he only dealt with the relative irregularity. 
Our approach is more natural, as the unitary flat subbundle in the second Fujita decomposition sits naturally in the Harder-Narasimhan filtration of the Hodge bundle: see \Cref{rem: ok}.
Moreover, to obtain the strong inequalities results in \Cref{teo: principale}, we use crucially the new Clifford-type inequality for canonical curves proved by the author  together with Riva in \cite{R-S}.

With Xiao's approach it is very hard to prove something on the extremal cases: for example the fact that surfaces attaining the classical slope \eqref{eq: canslope} are hyperelliptic  cannot be derived from the Xiao method. 
In  \Cref{teo: principale}, instead,  we give much information on the cases reaching the equalities.
This is done using  other techniques: in particular some results in \cite{GTS},  the algebraic index theorem and results of Lu and Zuo \cite{L-Z-hyper} and Konno \cite{konno3}.

\medskip

\noindent{\bf Acknowledgements} This article would not have been written without the encouragement and stimulation of Gian Pietro Pirola, whom I thank heartily. I thank the referees for their careful reading of the paper and for many useful suggestions, in particular I owe to one of them a sharpening of \Cref{teo: principale}. I am partially supported by the PRIN project 20228JRCYB ``Moduli spaces and special varieties'' and by GNSAGA - INdAM.

\section{Fibred surfaces and their geography}\label{sec: preliminaries}
We call a \textit{fibred surface} or sometimes simply a {\em fibration} the data of a  morphism $f\colon S\rightarrow B$ from a smooth projective surface $S$ to a smooth projective curve $B$ which is proper, surjective and with connected fibres.
We denote by $b=g(B)$  the genus of the base curve. A general fibre $F$ is a smooth curve and its genus $g=g(F)$ is by definition the genus of the fibration. From now on, we consider fibrations of genus $g\geq 2$.

We say that $f$ is \textit{relatively minimal} if it does not contain any $(-1)$-curves in its fibres. 
This condition is equivalent to the relative canonical divisor $K_{f}=K_S-f^*K_B$ being relatively nef. 
Given any fibration of genus $g\geq 2$ we can contract all the $(-1)$-curves contained in the fibres and obtain a unique relatively minimal genus $g$ fibration $\overline f\colon \overline S\to B$. 

Given a fibred surface, it is natural to define some {\em relative numerical invariants}, as follows:
\begin{itemize}
\item $K_f^2=K_S^2-8(g-1)(b-1)$ the self-intersection of the relative canonical divisor;
\item $\chi_{f}:=\chi(\mathcal{O}_S)-(g-1)(b-1)$ the relative Euler characteristic;
\item $e_{f}:= e(S)-e(B)e(F)=e(S)-4(g-1)(b-1)$ the relative topological characteristic (with $e(S)=c_2(S)$ the topological characteristic of $S$);
\item $q_{f}:=q-b$ the relative irregularity, with $q=h^1(S,\mathcal{O}_{S})$ irregularity of $S$ (see also \Cref{rem: ossqf});
\item the gonality  $gon(F)$ of the general fibre and its Clifford index $c_f=\cliff(F)$.
\end{itemize}
From the Groethendieck-Riemann-Roch theorem we have the relative Noether equality:
\begin{equation}\label{eq: noet}
12\chi_f = K_f^2 +e_f .
\end{equation}

We say that a fibred surface is \textit{smooth}  if every fibre is smooth;
 \textit{isotrivial} if all smooth fibres are mutually isomorphic;
 \textit{locally trivial} if $f$ is smooth and isotrivial (equivalently  $f$ is a fibre bundle by the Grauert-Fisher Theorem);
 {\em trivial} if $S$ is birationally equivalent to $F\times B$ and $f$ corresponds to the projection on $B$. 
 If $f$ is relatively minimal this is equivalent to $S=F\times B$.
Eventually, we say that $f$ is a {\em Kodaira fibration} if it is smooth and non-isotrivial.

Let us summarize what we know about the non-negativity of these invariants.
\begin{theorem}\label{teoremone}
Let $f\colon S\to B$  be a  fibred surface. The following results hold:
\begin{enumerate}
\item[(1)] Suppose that $f$ is relatively minimal. Then  $K_{f}^2\geq 0$ and $K_{f}^2=0$ if and only if $f$ is locally trivial;
\item[(2)]\label{ii} $\chi_{f}\geq 0$ and $\chi_f =0$ if and only if $f$ is locally trivial;
\item[(3)] \label{ef} $e_{f}\geq 0$ and $e_{f}=0$ if and only if $f$ is smooth;
\item[(4)]\label{beau} $q_f\leq g$ and the equality holds if and only if $f$ is trivial (see \cite{Beau} and also  \Cref{rem: loc}).
\end{enumerate}
\end{theorem}

\begin{definition}
The rank $g$ vector bundle $f_{*}\omega_{f}$ is called the \textit{Hodge bundle} of the fibred surface.
\end{definition} 

By using the Grothendieck-Riemann-Roch theorem or  Leray's spectral sequence, we see that 
$\deg f_{*}\omega_{f}=\chi_f$
for any fibration $f$.

\begin{definition}\label{def: nef}
A vector bundle $\cE$ over a smooth curve is called {\em nef} if  one of the following equivalent properties holds:
(1)  the corresponding tautological sheaf $\cO_{\pr(\cE)}(1)$ is nef;
(2) all its quotients have non-negative degree. 
\end{definition}

Let us recall some important results on the Hodge bundle.

\begin{theorem}\label{teoremone2}
The Hodge bundle of a fibration $f\colon S\to B$ can be decomposed in two ways as a direct summand of vector sub-bundles as follows:
\begin{itemize}
\item (First Fujita decomposition \cite{Fuj1})
\begin{equation}\label{first fujita}
f_{*}\omega_{f}=\mathcal{O}_{B}^{\oplus q_{f}}\oplus \mathcal{G}, 
\end{equation}
where $\mathcal{G}$ is nef and  $H^0(B,\mathcal{G}^\vee)=H^1(B, \cG\otimes \omega_B)=0$;
\item  (Second Fujita decomposition \cite{Fuj2} \cite{CD2})
\begin{equation}\label{second fujita} 
f_{*}\omega_{f}=\mathcal{A}\oplus \mathcal{U},
\end{equation}
where $\mathcal{A}$ is ample and $\mathcal{U}$ is unitary flat.
\end{itemize}
\end{theorem}
\begin{remark}\label{rem: ossqf}
We see from the above results that $q_f$ is the rank of the biggest trivial subbundle of $\cE$. Moreover, we see that $\chi_f=\deg \cG=\deg \cA$. 
We can also derive directly from the first decomposition that $q_f\leq g$ and that if $q_f=g$ then  $\chi_f=0$, and so $f$ is locally trivial (but Beauville's result (4) of \Cref{teoremone} is stronger).
\end{remark}
Following \cite{GTS}, we define the \textit{unitary rank} $u_{f}$ of the fibred surface to be the rank of $\cU$. 
\begin{remark}\label{rem: loc}
Comparing the two decompositions, since every trivial bundle is unitary flat, we have that:
$$\mathcal{O}_{B}^{\oplus q_{f}}\subseteq \mathcal{U}, $$
and then it holds that $q_{f}\leq u_{f}$. Observe that $\cU$ is the biggest degree $0$ sub-bundle of $f_*\omega_f$, so $u_{f}=g$ if and only if $\chi_f=0$ (equivalently $f$ is locally trivial): compare with (2) of  \Cref{teoremone}.
\end{remark}

One of the most important areas of study in this field is the so-called {\em geography} of fibred surfaces: the study of the inequalities holding between the invariants. We will list three geographical inequalities. For very nice expositions of these problems see \cite{MLP} and \cite{BCP}.

\subsection{Inequalities between the invariants: bounds on the slope $K_f^2/\chi_f$}

\begin{definition}
Let $f\colon S\to B$ be a locally non-trivial fibred surface. The {\em slope of $f$} is  the ratio $K_f^2/\chi_f$.
\end{definition}
Usually in the literature the   {\em slope inequalities} are the lower bounds on the slope. 

We start by recalling two upper bounds. 
First of all, from Noether's relation \eqref{eq: noet} and \eqref{ef} of  \Cref{teoremone} we have:
\begin{theorem}
For any fibred surface of genus $g\geq 2$, we have 
\begin{equation}\label{eq: 12}
K_f^2\leq 12\chi_f.
\end{equation}
Moreover, suppose that  $f$ is relatively minimal and  locally non-trivial. 
Then  the equality in \eqref{eq: 12} holds if and only if  $f$ is a Kodaira fibration.
\end{theorem}
\begin{remark}
Given any relatively minimal locally non-trivial fibred surface $f\colon S \to B$ over a base of genus $b\geq 2$, the surface $S$ is minimal and of general type. So, the Bogomolov-Miyaoka-Yau inequality holds: $K_S^2\leq 9\chi(\cO_S)$. This implies that 
\[K_f^2\leq 9\chi_f+(g-1)(b-1).\]
\end{remark}
\begin{remark}
The ratio 8 is particularly significant for the geometry of the surface. Indeed, consider the topological index $\tau(S)$  of the surface, i.e., the index of the intersection form on $H^2(S,\mathbb{C})$ 
(the difference between the number of its positive and negative eigenvalues).
We have that 
 $\tau(S) = \frac{1}{3} (K_S^2 -2c_2(S)) $.
 Therefore, using  the equalities $K_f^2=K_S^2-8(g-1)(b-1)$ and $\chi_{f}:=\chi(\mathcal{O}_S)-(g-1)(b-1)$ and the Noether equality $12\chi(\cO_s)=K_S^2+c_2(S)$, we have:
\[K_f^2-8\chi_f=K_S^2-8\chi(\cO_S)=\frac{K_S^2-2c_2(S)}{3}=\tau(S).\]
\end{remark}

Concerning  the lower bounds on the slope, the first and most famous result is the so-called classical slope inequality:
\begin{theorem}\label{teo: slope}[Xiao \cite{Xiao}, Cornalba-Harris \cite{C-H, Sto}, see also \cite[Thm. 4.5]{ACGH2}, Konno \cite{konno3}]

\noindent Let $f\colon S\to $ be a relatively minimal fibred surface, then the following inequality holds:
\begin{equation}\label{eq: canslope}
K_f^2\geq \frac{4g-4}{g}\chi_f.
\end{equation}
If $f$ is locally non-trivial and  the equality in \eqref{eq: canslope} is reached, then $F$ is hyperelliptic and $q_f=0$.
\end{theorem}

\begin{remark}
Even from the classical results of  \Cref{teo: slope}, it is natural to ask for lower  bounds on the slope  increasing with $c_f$ and with $q_f$. Many authors have contributed to this topic; here we just recall \cite{konno-cliff}, \cite{SB},  \cite{BS-trigonal}, \cite{LZ3}.
\end{remark}

\subsection{Inequalities between the invariants: Xiao's type of bounds}
As mentioned before, we have that  $q_f\leq g$, and the equality holds if and only if $f$ is trivial. 
Therefore, it is natural to ask: if $f$ is not trivial, what is the inequality between $g$ and $q_f$?
From now on we assume $f$ is non-trivial. 
Xiao initiated the study in \cite{xiao-irregular} by proving that if $b=0$ 
\begin{equation}\label{eq: xiao bound}
q_f\leq \frac{g+1}{2}. 
\end{equation}
For arbitrary $b\geq 0$ he proved the bound 
\begin{equation}\label{eq: xiaobrutto}
q_f\leq \frac{5g+1}{6}.
\end{equation}
He also conjectured  that the bound  \eqref{eq: xiao bound} holds for any non-trivial fibration. 
Serrano \cite{serrano} proved that if  $f$ is isotrivial  (but not trivial), then \eqref{eq: xiao bound} holds.
Cai in  \cite{Cai} proved that if $f$ is non-isotrivial and the general fibre is either hyperelliptic or bielliptic, the same bound holds. 
In \cite{konno}, Konno improved the bound \eqref{eq: xiaobrutto} to: 
\begin{equation}\label{eq: konno}
q_f \leq g \frac{5g-2}{6g-3}.
\end{equation}
Pirola gave a genus $4$  counterexample to Xiao's conjecture in \cite{pirola}. 
More counterexamples have been found by Albano and Pirola in \cite{AP}.
In  \cite{BGN}  Barja, Gonz\'alez and Naranjo proved the following: if $f$ is non-isotrivial, then $q_f\leq g-c_f$.
They made the following modified conjecture:
\begin{conjecture}[Barja, Gonz\'alez, Naranjo; Modified Xiao's conjecture for the relative irregularity]\label{conj:mxiao}
For any non-isotrivial fibred surface $f\colon S\to B$ of genus $g\geq 2$ it holds 
\begin{equation} \label{eq:BGN}
q_f\leq \left\lceil \frac{g+1}{2} \right\rceil.
\end{equation}
\end{conjecture}
Favale-Naranjo-Pirola \cite{FNP} proved the stronger inequality $q_f\leq g-c_f-1$ for families of plane curves of degree $d \geq 5$.
It is natural to ask whether these bounds can be extended to the unitary rank $u_f$. 
In the paper \cite{GTS} the author with Gonz\'alez and Torelli proved that :
$(i)$  $u_f \leq g - c_f;$
$(ii)$ if the general fibre is a plane curve
of degree $d \geq 5$, then $ u_f \leq  g - c_f - 1.$
\begin{remark}\label{rem: confronto}
If $\cU$ has finite monodromy, then up to an \`etale base change it becomes trivial, and so coincides with $q_f$. However, the unitary summand $\cU$ of the Hodge bundle can have infinite monodromy as proved by Catanese and Dettweiler in \cite{CD2}. Moreover, the modified Xiao's conjecture does  {\em not} hold for $u_f$:  this follows from \cite[Remark 38]{CD2} and from the construction of Lu \cite{Lu}. See also \cite[Conjecture 1.1]{GTS}.
 \end{remark}


\section{Main technical results}\label{konno-method}

\subsection{The Harder-Narasimhan filtration of a vector bundle}
The slope of a vector bundle $\cE$ over a curve $B$  is defined as $\mu(\cE):=\deg(\cE)/\rk(\cE)$.
The vector bundle $\cE$ is called {\em semistable (resp. stable)} if for any vector sub-bundle $\cF\subsetneq \cE$ we have $\mu(\cF)\leq\mu(\cE)$ (resp. $\mu(\cF)<\mu(\cE)$).

Let us briefly recall the notion of the Harder-Narasimhan filtration of  a vector bundle $\cE$ over a curve $B$  \cite{H-N}. 
It is the unique filtration of subbundles
$$0=\cE_0\subsetneq \cE_1\subsetneq \ldots \subsetneq \cE_l=\cE$$
satisfying the following assumptions:
\begin{itemize}
\item $\cE_i/\cE_{i-1}$ is  a $\mu$-semistable vector bundle for any $i=0, \ldots l$,
\item if we set $\mu_i:= \mu(\cE_i/\cE_{i-1})$  for any $i=1, \ldots l$ we have that $\mu_i<\mu_{i-1}$.
\end{itemize}
Note that $\mu_1>\mu(\cE)>\mu_l$, unless $\cE$ is $\mu$-semistable, in which case $1=l$ and $\mu(\cE)=\mu_1$. 

Let us set $r_i:=\rk(\cE_i)$. 
Setting $r_0:=0$, we can express the degree of $\cE$ as a combinations of the $\mu_i$'s and the $r_i$'s: 
\begin{equation}\label{eq: grado}
\deg(\cE)=\sum_{i=1}^{l}\mu_i(r_{i}-r_{i-1}).
\end{equation}
Let us call $H$  a tautological divisor and $\Sigma$ a fibre of the projective bundle $\pi\colon \pr_B(\cE)\to B$.
Recall that $\pic(\pr_B(\cE))=\pi^*\pic(B)\oplus \Z[H]$. 
The following result of Miyaoka-Nakayama \cite{Miy, nak} determines the form of the nef and pseudo-effective cones of $\pr_B(\cE)$.
\begin{theorem}\label{thm: M-N}[Miyaoka-Nakayama]\\
With the above notations, given $k\in \Z$ and  $D$ a divisor on $B$,
\begin{enumerate}
\item The  divisor  $kH-\pi^*D$ is pseudo-effective (i.e., it is a limit of effective divisors) if and only if $k\geq 0$ and $\deg D\leq k\mu_1$.
\item The  divisor  $kH-\pi^*D$ is nef if and only if $k\geq 0$ and $\deg D\leq k\mu_l $.
\end{enumerate}
\end{theorem}
\begin{remark}
Although the theorem is stated for divisors, it can be readily extended to $\Q$-divisors, and it is  this version that we will use.
\end{remark}

\subsection{The Harder-Narasimhan filtration of the Hodge bundle}\label{ssec: HN}
Let us consider the case of a fibred surface $f\colon S\to B$, and consider its Hodge bundle  $\cE=f_*\omega_f$.
As $\cE$ is nef by  \Cref{teoremone2}, we have that $\mu_l$ is greater than or equal to $0$. 
\begin{remark}\label{rem: ok}
The unitary flat summand $\cU$ in \eqref{second fujita} is the biggest degree $0$ subbundle of $\cE$. 
So, in particular we have that $\mu_{l}=0$ if and only if $\mathcal U\not=0$ if and only if $u_f\not =0$.
Moreover, if $\mu_l=0$, then $\cE_{l-1}$ is precisely $\mathcal A$, if $\mu_l>0$, then $\mathcal A$ is the whole Hodge bundle.
\end{remark}
\begin{remark}\label{rem: vabbuo}
Observe moreover that we have in general:
\[\mu_l\leq \frac{\chi_f}{g-u_f} \leq  \frac{\chi_f}{g-q_f} \leq \frac{\chi_f}{g}\leq \mu_1\]
and in particular:
\begin{itemize}
\item[(i)] $\mu_1=\chi_f/g=\mu(f_*\omega_f)$ if and only if $f_*\omega_f=\cA$ is semistable.
\item[(ii)] $\mu_1=\chi_f/(g-u_f)=\mu(\cA)$ if and only if $\cA$ is semistable.
\item[(iii)] $\mu_1=\chi_f/(g-q_f)=\mu(\cG)$ if and only if $u_f=q_f$ and $\cA$ is semistable.
\end{itemize}
\end{remark}

\subsection{Xiao's set-up}\label{ssec: xiao}

We follow the  set-up given by Xiao in his seminal paper  \cite{Xiao}. 
Let us consider the Hodge bundle $\cE=f_*\omega_f$, and its Harder-Narasimhan filtration as in  \Cref{ssec: HN}. 
For any $i=1,\ldots l$, the sheaf homomorphism 
\[f^*\cE_i\longrightarrow f^*f_*\omega_f\longrightarrow \omega_f\]
induces a rational map $\phi_i\colon S\dasharrow \pr_B(\cE_i)$. 
Let $H_i$ be the tautological divisor on $\pr_B(\cE_i)$ and let $M_i:=\phi_i^*(H_i)$ on $S$.
\begin{remark}\label{rem: M-N}
By  \Cref{thm: M-N}, observing that $\mu_i$ is the final slope of $\cE_i$ for all $i$, we have that $H_i-\mu_i\Sigma_i$ is a nef $\Q$-divisor and  $H_i-\mu_1\Sigma_i$ is a pseudoeffective $\Q$-divisor on $\pr_B(\cE_i)$.
\end{remark}
Hence, $\phi_i^*(H_i-\mu_i\Sigma_i)= M_i-\mu_i F$ is a nef $\Q$-divisor on $S$ and $\phi_i^*(H_i-\mu_1\Sigma_i)= M_i-\mu_1F$ is a pseudoeffective $\Q$-divisor on $S$.
Moreover, for any $i=1,\ldots l$, call $Z_i$ the divisorial base locus of $\phi_i$, i.e., the effective divisor in $S$ such that  
$\cE_i\subseteq  f_*\omega_f(-Z_i)$ and such that 
the evaluation homomorphism
\[f^*\cE_i\longrightarrow f^*f_*\omega_f(-Z_i)\longrightarrow \omega_f(-Z_i)\]
is surjective in codimension $2$. 
Of course
\[Z_1\geq Z_2\geq \ldots \geq Z_l=0.\]
Due to the importance of this result in what follows, we give the proof of this lemma:
\begin{lemma}\label{lem: punto}
With the above notations, for any $i=1,\ldots l$
the divisor $ K_f $ is numerically equivalent to $ M_i+ Z_i$, 
\end{lemma}
\begin{proof}
The sheaf homomorphism $f^*\cE_i\to f^*f_*\omega_f(-Z_i)\to \omega_f(-Z_i)$ is surjective in codimension 2. 
Let us consider a chain of blow ups  $\rho\colon \widetilde S\to S$  resolving the points of indeterminacy of the associated map $\phi_i$: there exists a morphism $\tilde \phi_i\colon \widetilde S\to \pr_B(\cE_i)$ such that $ \tilde \phi_i=\phi_i\circ \rho$.
We thus have the following commutative diagram:
\[
\xymatrix{
\widetilde S\ar_\rho[d]\ar[rrd]^{\widetilde{\phi_i}}&&\\
S\ar_f[rd]\ar^{\phi_i}@{-->}[rr]&&\pr_B(\cE_i)\ar^\pi[dl]\\
&B&}
\]
Let $\tilde f=f\circ\rho\colon \widetilde S\to B$ be the associated fibration;  note that it is not relatively minimal unless $\phi_i$ already is a morphism. 
Let $E_i$ be the exceptional divisor of $\rho$. 

The morphism $\tilde \phi_i$ is associated to the line bundle $\mathcal{L}=\rho^*\omega_f(-Z_i)\otimes\cO_{\widetilde S}(-E_i)$ on $\widetilde S$, and we have the surjective sheaf homomorphism:
\[\tilde f^*\cE_i\longrightarrow {\tilde f}^*{\tilde f}_*\mathcal{L}\longrightarrow \mathcal{L} .\]
So 
we have that $\mathcal{L}\cong \widetilde{\phi_i}^*(\cO_{\pr_B(\cE_i)}(H_i))$.
Now, given any irreducible curve $C\subset S$, we have by using the push-pull formula that 
\begin{equation*}
\begin{split}
M_i\cdot C&=\phi_i^*(H_i)\cdot C=\tilde \phi_i^*(H_i)\cdot \rho^*(C)=\rho^*(K_f-Z_i)\cdot \rho^*(C)+E_i\cdot \rho^*(C)=\\
&= \rho^*(K_f-Z_i)\cdot \rho^*(C)=(K_f-Z_i)\cdot C,\\
\end{split}
\end{equation*}
as wanted.
\end{proof}
 As a direct consequence of the previous discussion and of \Cref{lem: punto} we have the following:
 \begin{proposition}\label{prop: xiao}[Xiao]
 With the above notations, we have that for any $i=1,\ldots l$
 \begin{itemize} 
 \item[(i)]  $K_f-Z_i- \mu_iF$ is nef;
 \item[(ii)] $K_f-Z_i-\mu_1F$ is pseudoeffective.
 \end{itemize}
 \end{proposition}
 \begin{remark}\label{rem: nofibres}
 Note that the divisors $Z_i$'s can  have vertical components. In the case when they have fibres as components, say $Z_i=bF+Z'_i$ with $b\geq 0$ and the support of  $Z_i'$ containing no fibre, then we can observe the following:
 \begin{itemize} 
 \item[$(i')$]  $K_f-Z'_i- \mu_iF$ is nef;
 \item[($ii'$)] $K_f-Z'_i-\mu_1F$ is pseudoeffective.
 \end{itemize}
 Indeed, for $(i')$,
 let $C$ be any irreducible curve in $S$. If $C$ is vertical then $C\cdot F=0$ and so 
 \[(K_f-Z'_i- \mu_iF)\cdot C=(K_f-Z_i- \mu_iF)\cdot C\geq 0.\] In the case  where$C$ is horizontal then $C\cdot F>0$, so 
 \[(K_f-Z'_i- \mu_iF)\cdot C=(K_f-Z_i- \mu_iF)\cdot C+bF\cdot C\geq 0.\]
 For $(ii')$ just observe that $K_f-Z'_i-\mu_1F= K_f-Z_i-\mu_1F+bF$, so it is certainly pseudoeffective.
 For this reason in the following we can assume that the $Z_i$'s can only have vertical components that are not fibres.
 We cannot suppose that $Z_i$ has no vertical components at all, because $(i')$ is not guaranteed if we eliminate those components.
 
 \smallskip
 
 Another approach to get rid of fibres in the base divisors --the one used by Xiao in \cite{Xiao}-- is to tensor $\cE_i$ with a sufficiently ample line bundle $\cB$ on $B$ and to consider $L$ the linear subsystem of $|\omega_f\otimes f^*\cB|$ corresponding to sections in $H^0(B, \cE_i\otimes \cB)$. Then one defines  $Z_i$ to be the fixed part of $L$ and verifies that this construction is independent on the choice of the ample line bundle $\cB$.
 \end{remark}
\begin{notation}
Let us call $d_i=M_i\cdot F$, the degree of the linear series $|V_i|$ induced by $M_i$ on the general fibre, and let $a_i=(Z_i\cdot F)=2g-2-d_i$. The divisor ${Z_i}_{|F}$  on $F$ is the base locus of the linear series $|V_i|$ seen as a subcanonical system.

Recall that $r_i:=\rk \cE_i$, for $i=1,\ldots, l$.  The linear series $|V_i|$ is associated to a base-point-free subcanonical $g^{d_i}_{r_i-1}$ on $F$.
\end{notation}
\begin{remark}\label{rem: piccolo ma}
Observe that $d_1\leq d_2\leq \ldots \leq d_l=2g-2$. Also note that $d_i>0$ (i.e., $a_i<2g-2$) if $i>1$; indeed, if $i>1$ the linear series $|V_i|$ has dimension $r_i-1>0$. 
Moreover, $d_1=0$ if and only if $r_1=1$.
\end{remark}


\subsection{A new nef divisor on $S$}
Till now we have been working in the framework of Xiao's technique. 
Now we turn to an original idea of Konno.
%
\begin{proposition}[Konno, \cite{konno}]\label{prop: fond}
With the above notation,
let $Z_i=\sum_j m_jG_j$ be the irreducible decomposition of $Z_i$, and let $\alpha_i:=\max_j\{m_j\mid  F\cdot G_j>0\}$.
Then the $\Q$-divisor $tK_f+M_i-\mu_iF+Z_i$ is nef for any $t\geq \alpha_i$.
In particular, $(\alpha_i+1)K_f -\mu_iF\equiv \alpha_i K_f +M_i-\mu_iF+Z_i$ is a nef divisor.
\end{proposition}
\begin{proof}
For the readers' convenience we elaborate here Konno's proof of  \cite[Lemma 2.1]{konno}. Fix $i\in \{1,\ldots , l\}$. Let $Z_i=\sum_j m_jG_j$ be the irreducible decomposition of $Z_i$, and let 
\[\alpha_i:=\max_j\{m_j\mid  F\cdot G_j>0\}.\] 
In other words, $\alpha_i$ is the maximal multiplicity of the non-vertical components of $Z_i$. 
From  \Cref{rem: nofibres} we can assume that $Z_i$ contains no fibres, but it can contain some components of fibres. 
From  \Cref{teoremone},  \Cref{lem: punto} and \Cref{prop: xiao} we know  that $K_f$ and $M_i-\mu_iF$  are nef. 
Let us now consider any irreducible curve $C\subset S$. If $C$ is not a component of $Z_i$, then $C\cdot Z_i\geq 0$ and so $(\alpha_i K_f +M_i-\mu_iF+Z_i)\cdot C\geq 0$. 
Let us now assume that $C=G_j$ for some $j$. If $C$ is vertical, $C\cdot F=0$ and so 
$$(\alpha_i K_f +M_i-\mu_iF+Z_i)\cdot C=((\alpha_i+1)K_f -\mu_iF)\cdot C=(\alpha_i+1)K_f\cdot C\geq 0.$$
If on the contrary, $C\cdot F>0$, then the restriction of $f$ to $C$ is a finite surjective morphism $\pi=f_{|C}\colon C\to B$ of degree $d=C\cdot F$. We have that 
$$(K_f+C)\cdot C=(K_S+C-f^*K_B)\cdot C=\deg \omega_C - d(2b-2)=\deg r(\pi)\geq 0,$$
by the Riemann-Hurwitz formula, where $r(\pi) $ denotes the  ramification divisor of $\pi$ on $C$. 
Now we can rearrange the divisor $\alpha_i K_f +Z_i$ as follows:
$$ \alpha_i K_f +Z_i=(\alpha_i-m_j)K_f+m_j(K_f+C)+ (Z_i-m_jC);$$
and so we have 
$$(\alpha_i K_f +M_i-\mu_iF+Z_i)\cdot C=(\alpha_i-m_j)K_f\cdot C+m_j( K_f+C)\cdot C+  (Z_i-m_jC)\cdot C +(M_i-\mu_iF)\cdot C.$$
All the contributions are greater than or equal to zero, by the following observations:
\begin{itemize}
\item[(i)]$\alpha_i-m_j\geq 0$ by definition of $\alpha_i$, and $K_f$ is nef. 
\item[(ii)]$m_j( K_f+C)\cdot C=m_j\deg r(\pi)\geq 0$.
\item[(iii)]$Z_i-m_jC $ is an effective divisor which does not contain $C$ in its support, so it has non-negative intersection with $C$.
\item[(iv)] $M_i-\mu_iF$ is nef.
\end{itemize}
The proposition is thus concluded.
\end{proof}
\begin{remark}\label{rem: novita'}
It is important to remark that this idea of Konno introduces a real novelty with respect to  Xiao's method; indeed, the latter method constructs nef divisors on $S$ by  pulling back  (extremal) nef divisors on $\pr_{B}(\cE_i)$, so it obtains  divisors of the form $M_i-\mu_iF\equiv K_f -Z_i-\mu_iF$. In particular  it does not give information on the nefness of divisors of the form $K_f-bF$, unless $Z_i=0$.
\end{remark}
\begin{remark}\label{rem: alpha}
Notice that $\alpha_i \leq a_i$, and if $\alpha_i=a_i$, then necessarily $Z_i=a_i\Gamma +Z'$, where $\Gamma$ is a section of $f$ and $Z'$ is an effective vertical divisor, and moreover $a_i\Gamma $ is contained in all the $Z_j$ with $j\leq i$.
\end{remark}

\begin{lemma}\label{lem: chiave}
With the above notations, we have that for any $i=1,\ldots l$  
\begin{equation}\label{eq: chiave}
K_f\cdot Z_i\geq\mu_i \frac{a_i}{\alpha_i+1}.
\end{equation}
\end{lemma}
\begin{proof}
We consider the rational map $\phi_i\colon S\dasharrow \pr_B(\cE_i)$ defined at the beginning of  \Cref{ssec: xiao}.
As proved above, $K_f$ is numerically equivalent to $M_i+Z_i$, with $Z_i$ effective and $M_i-\mu_iF$ nef. 
By  \Cref{prop: fond} we have that $(a_i+1)K_f-\mu_iF\equiv a_iK_f+M_i+Z_i-\mu_iF$ is a nef divisor. As $a_i\geq \alpha_i$, also 
 $(\alpha_i+1)K_f-\mu_iF\equiv \alpha_iK_f+M_i+Z_i-\mu_iF$ is nef.
We thus have 
$((\alpha_i+1)K_f-\mu_iF)\cdot Z_i\geq 0$, and ultimately  $K_f\cdot Z_i\geq a_i \mu_i/(\alpha_i+1)$. 
\end{proof}
We now  derive a general slope inequality:
\begin{lemma}\label{lem: konnoplus}
Let $f\colon S\to B$ be a relatively minimal locally non-trivial fibred surface of genus $g\geq 2$. Let $a:=a_1$, and $\alpha:=\alpha_1$.
The following inequality holds:
\begin{equation}\label{eq: 2}
K_f^2\geq \left(\frac{(\alpha+1)(4g-4-a)+a}{\alpha+1}\right)\mu_1.
\end{equation}
\end{lemma}
\begin{proof}
We just consider the first rational map $\phi_1\colon S\dasharrow \pr_B(\cE_1)$. By  \Cref{lem: chiave} we have that $K_f\cdot Z_1\geq a\mu_1/(\alpha+1)$.
Now observe that:
\begin{equation}\label{eq: principal}
\begin{split}
K_f^2&=K_f\cdot (M_1+Z_1)=K_f\cdot M_1+K_f\cdot Z_1=K_f\cdot (M_1-\mu_1F)+K_f\cdot \mu_1F+K_f\cdot Z_1=\\
&=K_f-\mu_1F\cdot (M_1-\mu_1F)+\mu_1F\cdot (M-\mu_1F)+(K_f\cdot \mu_1F+K_f)\cdot Z_1.
\end{split}
\end{equation}
Now, we observe that $K_f-\mu_1F\equiv M_1-\mu_1F+Z_1$ is a pseudo-effective divisor, because $M_1-\mu_1F$ is pseudo-effective by  \Cref{thm: M-N} (as $\mu_1\geq \mu_l$), and $Z_1$ is effective by its very definition. So, the first term $(K_f-\mu_1F)\cdot (M_1-\mu_1F)$ is greater than or equal to $0$, and we obtain 
\begin{equation}\label{eq: int}
\begin{split}
K_f^2&\geq \mu_1(F\cdot M_1)+\mu_1(K_f\cdot F) + (K_f\cdot Z_1)\geq\quad \quad \quad \quad \quad \\
&\geq \mu_1\left((2g-2-a+2g-2)+\frac{a}{\alpha+1}\right)=\mu_1\left(\frac{(\alpha+1)(4g-4-a)+a}{\alpha+1}\right),
\end{split}
\end{equation}
as claimed.
\end{proof}

\begin{remark}\label{rem: generale}
More generally one can prove that for any $i=1,\ldots ,l$ 
\begin{equation}\label{eq: 3}
K_f^2\geq \left(\frac{(\alpha_i+1)(4g-4-a_i)+a_i}{\alpha_i+1}\right)\mu_i.
\end{equation}
\end{remark}

\section{Some new slope inequalities}

We prove in this section several slope inequalities involving the unitary rank $u_f$, and some also involving $c_f$, and the rank $r_1$ of the first subbundle $\cE_1$ of the Harder-Narasimhan filtration of the Hodge bundle.
\subsection{First applications}
Even avoiding the elaboration made in Lemmas \ref{lem: chiave} and \ref{lem: punto},
one can derive  directly from  \Cref{prop: fond} some interesting slope inequalities, 
which will be very important for the applications in the following sections. 

\begin{theorem}\label{teo: konnoplus}
Let $f\colon S\to B$ be a relatively minimal locally non-trivial fibred surface of genus $g\geq 2$. Let $\alpha:=\alpha_1$ be the maximal multiplicity of the non-vertical irreducible components of $Z_1$.
The following inequality holds:
\begin{equation}\label{eq: chiave}
K_f^2\geq \left(\frac{2g-2}{\alpha+1}+2g-2\right)\mu_1.
\end{equation}
 In particular, the following inequality holds:
\begin{equation} \label{eq: konnoplus}
K_f^2 \geq \frac{4g(g-1)}{(2g-1)(g-u_f)}\chi_f.
\end{equation}
Moreover, the equality holds in \eqref{eq: konnoplus} if and only if $u_f=g-1$ and $K_f $ is numerically equivalent to  $(2g-2)\Gamma+ \chi_f F=(2g-2)\Gamma+ \mu_1 F$, where $\Gamma $ is a section of $f$.
\end{theorem}
\begin{proof} 
First observe that the assumption that $f$ is locally non-trivial implies that $\chi_f>0$ and so equivalently (\Cref{rem: ok}) $u_f\leq g-1$. 
From \Cref{prop: fond} we have that
$(\alpha+1)K_f -\mu_1F$ is a nef divisor on $S$, and by \Cref{prop: xiao} $K_f-\mu_1F\equiv M_1-\mu_1F+Z_1$ is a pseudo-effective divisor, so their intersection is greater than or equal to $0$:
 \[((\alpha+1)K_f -\mu_1F)\cdot ( K_f-\mu_1F)=(\alpha+1)K_f^2-(2g-2)(\alpha+2)\mu_1\geq 0.\]
 So 
\[K_f^2\geq (2g-2)\frac{\alpha+2}{\alpha+1}\mu_1.\]
Therefore, the inequality \eqref{eq: chiave} is proved.
Now just observe that the function 
\[t\mapsto \frac{t+2}{t+1}\] 
is strictly decreasing and $\alpha\leq 2g-2$ so we obtain
\[K_f^2\geq \frac{(2g-2)(2g)}{2g-1}\mu_1\geq \frac{(2g-2)(2g)}{2g-1}\frac{\chi_f}{g-u_f},\]
where the last inequality follows from \Cref{rem: vabbuo}.

Suppose now that the equality holds in \eqref{eq: konnoplus}.  
Then the equalities also hold in all the passages of the proof. 
In particular:
\begin{itemize}
\item $\alpha =a=2g-2$, i.e $d_1=0$, so by \Cref{rem: M-N} we have that  $r_1=1$; 
\item $\mu_1=\chi_f/(g-u_f)=\mu (\cA)$, so by $(ii)$ of  \Cref{rem: vabbuo} $\cA$ is semistable.
\end{itemize}  
Putting all together, we have that necessarily the Harder-Narasimhan filtration has lenght $2$, that $g-u_f=\rk\cA=\rk\cE_1=1$ and that $\mu_1=\chi_f$. 

The fact that $\alpha=a=2g-2$ implies by \Cref{rem: alpha} that the divisor $Z_1$ 
is necessarily  $(2g-2)\Gamma+Z'$, where $\Gamma $ is a section of $f$, and $Z'$ is an effective  vertical divisor. 

Now we prove that if the equality holds in \eqref{eq: konnoplus} then $Z'=0$.

First of all observe that, via the same arguments used in  \Cref{rem: nofibres}, we see that $Z'$ does not contain any fibre in its support
 
We now prove that there are no vertical components. Suppose on the contrary $Z'\not=0$, and let $D$ be an irreducible component of $Z'$. 
We have the equality in \eqref{eq: konnoplus}
\[K_f^2 =\frac{4g(g-1)}{(2g-1)(g-u_f)}\chi_f;\]
On the other side $K_f^2=\frac{4g(g-1)}{2g-1}\mu_1+4(g-1)\Gamma \cdot Z'$, 
so we need to have $\Gamma \cdot Z'=0$ which implies $K_f\cdot Z'=Z'^2$. But $Z'^2\leq 0$ because $Z'$ is vertical, and so $K_f\cdot Z'=Z'^2=0$. 
The algebraic index theorem  \cite[Thm.~4.14]{ACGH2}, gives a contradiction: $K_f^2>0$ implies that $Z'^2<0$.
 
So, we eventually get that $K_f$ is necessarily numerically equivalent to $(2g-2)\Gamma+ \mu_1 F=(2g-2)\Gamma+ \chi_f F$. 

Now we prove the converse implication. Suppose that $u_f=g-1$ and  $K_f$ is  numerically equivalent to $(2g-2)\Gamma+ \chi_f F$, with $\Gamma$ a section of $f$. 
Then 
$K_f^2=2(2g-2)\chi_f+4(g-1)^2\Gamma^2$. Now, observe that $K_f\cdot \Gamma =-\Gamma^2$ by adjunction, and that  --by the algebraic index theorem again-- the intersection matrix between $K_f, \Gamma$ and $F$ has zero determinant:
\[
0=\det 
\begin{pmatrix}
K_f^2&K_f\cdot \Gamma &K_f\cdot F\\
K_f\cdot \Gamma & \Gamma^2& \Gamma\cdot F\\
K_f\cdot F &\Gamma\cdot F& F^2
\end{pmatrix}
=\det 
\begin{pmatrix}
K_f^2&-\Gamma^2 &2g-2\\
-\Gamma^2& \Gamma^2& 1\\
2g-2 &1& 0
\end{pmatrix}
=-2g(2g-2)\Gamma^2-K_f^2
\]
and hence $K_f^2=-4g(g-1)\Gamma^2$. Combining the two equalities we obtain precisely  
\[K_f^2=\frac{4g(g-1)}{2g-1}\chi_f.\]
\end{proof}

\begin{remark} 
The inequality \eqref{eq: konnoplus} is a  generalization of \cite[Lemma 2.7]{konno}.
For $u_f\ll g$ this result is  not very interesting; 
for example, in case $u_f=0$ \eqref{eq: konnoplus} gives as a coefficient  $(4(g-1))/(2g-1)$, which is strictly smaller than the one of the classical slope inequality \eqref{eq: canslope}. 
But for $u_f$ large with respect to $g$ the inequality  becomes interesting, as we will see in the next sections. 
\end{remark}

We now use again the algebraic index theorem in combination with Konno's method, to obtain a somewhat different inequality as follows. Let $\Sigma$ be an irreducible horizontal component of $Z_1$ of maximal multiplicity $\alpha=\alpha_1$. 
The restriction of $f$ to $\Sigma$ is a (possibly ramified) covering of $\pi\colon \Sigma\to B$, of degree say $\beta$.
\begin{theorem}
With the above notations,  the following inequality holds:
\begin{equation}\label{eq: hodge}
K_f^2\geq 4\frac{(4g-4-a)(g-1)(g-1+\beta)}{4(g-1)(g+1+\beta)-\beta^2\alpha}\mu_1
\end{equation}
\end{theorem}
\begin{proof}
First of all observe that $K_f\cdot \Sigma+\Sigma^2=\deg r (\pi)\geq0$, as proved in the argument of \Cref{prop: fond}. As $K_f^2>0$ we have that the intersection matrix between $K_f, \Sigma$ and $F$ has determinant greater than or equal to $0$:
\[
0\leq\det 
\begin{pmatrix}
K_f^2&K_f\cdot \Sigma &K_f\cdot F\\
K_f\cdot \Sigma & \Sigma^2& \Sigma\cdot F\\
K_f\cdot F &\Sigma\cdot F& F^2
\end{pmatrix}
=\det 
\begin{pmatrix}
K_f^2&K_f\cdot \Sigma &2g-2\\
K_f\cdot \Sigma& \Gamma^2& \beta \\
2g-2 &\beta& 0
\end{pmatrix}
\]
so we obtain 
\[0\leq 4\beta(g-1)K_f\cdot \Sigma-4(g-1)^2\Sigma^2-\beta^2K_f^2\leq (4\beta(g-1)+4(g-1)^2)K_f\cdot \Sigma -\beta^2K_f^2, \]
and so
\[K_f\cdot Z_1=K_f\cdot (\alpha \Sigma +Z')\geq \alpha K_f\cdot \Sigma \geq \frac{\alpha\beta^2}{4(g-1)(\beta+g-1)}K_f^2.\]
Now, plugging this inequality in \eqref{eq: int}, we obtain:
\[K_f^2\geq (4g-4-a)\mu_1+K_f\cdot Z_1\geq (4g-4-a)\mu_1+\frac{\alpha\beta^2}{4(g-1)(\beta +g-1)}K_f^2.\]
\end{proof}
\begin{remark}
Observe that in the case where $a=\alpha=2g-2$ (and thus $\beta=1$), we obtain \eqref{eq: konnoplus}. In general the inequality \eqref{eq: hodge} seems to be weaker than \eqref{eq: chiave}, unless one has specific informations on $\beta$ and $\alpha$.
\end{remark}
\begin{theorem}
\label{teo: troppo bella}
Let $f\colon S\to B$ be  a relatively minimal fibred surface of genus $g\geq 2$.
Assume that $Z_1$ doesn't have horizontal components, i.e. $a=\alpha=0$ (equivalently the linear system induced by the stalk of $\cE_1$ on the general fibre is base-point-free);
Then the following slope inequality holds:
\begin{equation}
\label{eq: troppo bella}
K_f^2\geq \frac{4(g-1)}{g-u_f}\chi_f.
\end{equation}
\end{theorem}
\begin{proof}
The inequality \eqref{eq: troppo bella} can  be deduced from \Cref{teo: konnoplus} setting $\alpha=0$ in \eqref{eq: konnoplus}. 
\end{proof}
\begin{remark}
In the case where $f$ is a hyperelliptic fibration, Lu and Zuo in \cite{L-Z-hyper} proved \eqref{eq: troppo bella} with the relative irregularity $q_f$ instead of the unitary rank $u_f$ In this case --again by Lu and Zuo's results,  in \cite[Theorem 4.7]{L-Z-hyper-oort}-- there is a finite \'etale base change such that the flat summand of the new Hodge bundle becomes  trivial, so  the inequality for $q_f$ implies the one with $u_f$.
\end{remark}
\begin{remark}
The inequality \eqref{eq: troppo bella} is false in general, as proved by the examples constructed by Lu in \cite{Lu}. Moreover, the assumption on the base locus of the system induced by $\cE_1$ is quite restrictive, and not very easy to verify in concrete cases. 
\end{remark}
\begin{remark}
The previous theorem is somewhat close to one due to the author with Barja in \cite[Remark 3.3]{SB}: there we proved that the  inequality \eqref{eq: troppo bella} holds if the ample summand $\cA$ of the Hodge bundle is semistable and its stalk induces a base-point-free system on the general fibres. Recall that $\cA$ is either the whole $f_*\omega_f$ or the one to last subsheaf of the Harder-Narasimhan filtration $\cE_l$. So, the assumptions of \cite[Remark 3.3]{SB} can be rephrased as follows: either the Hodge bundle is semistable (in this case by definition $Z_1=0$) or it has a Harder-Narasimhan filtration of length two, $\mu_2=0$ and $Z_1$ has no horizontal components. 
Therefore, the above Theorem implies \cite[Remark 3.3]{SB}.
\end{remark}

\subsection{Some new slope inequalities deriving from \Cref{lem: chiave}}
We now use  \Cref{lem: chiave}, and see that in many cases it gives stronger slope inequalities than the known ones.
\begin{theorem}\label{teo: tecnico}
Let $f\colon S\to B$ be a relatively minimal fibred surface of genus $g\geq 2$. 
Let us denote as usual $r_1=\rk\cE_1$.

The following inequality holds:
\begin{equation}\label{eq: r1}
K^2_f\geq \left(\frac{2g+2r_1-4}{g-u_f}+\frac{2g-2r_1}{(2g-2r_1+1)(g-u_f)}\right)\chi_f.
\end{equation}
Moreover, suppose that $f$ is locally non-trivial and that  the equality is reached. 
Then necessarily  $r_1=1$ and $f$ is a trigonal fibration, or $r_1=2$ and $g=2$.
In particular, for $g\geq 3$ and $r_1\geq 2$, inequality \eqref{eq: r1} is always strict.
\end{theorem}
\begin{proof}
 If $u_f>0$, then by  \Cref{rem: ok} we have $\cA=\cE_{l-1}$,
and  the Harder-Narasimhan filtration of $\cA$ is just the truncation of the one of $\cE$:
$$0=\cE_0\subset \cE_1\subset \ldots \subset \cE_{l-1}=\cA.$$
Let us recall  the inequality \eqref{eq: 2}:
\[K_f^2\geq \left(\frac{(\alpha+1)(4g-4-a)+a}{\alpha+1}\right)\mu_1,\]
where $a=a_1$ and $\alpha=\alpha_1$.
Observe that $a\geq \alpha$ and that the function 
\[f(t):=\frac{(t+1)(4g-4-t)+t}{t+1}\] 
is strictly decreasing.

If $r_1=1$, then $a=2g-2$, and we obtain \eqref{eq: konnoplus}, which coincides  with the  inequality \eqref{eq: r1}. 

If $r_1>1$, then by Clifford's theorem 
applied to the subcanonical linear system ${M_i}_{|F}$ we have $d_{1}\geq 2(r_{1}-1)$, hence $a_1= 2g-2-d_1\leq 2g-2r_1$.
We thus  have the following chain of inequalities:
\begin{equation*}
\begin{split}
K_f^2&\geq\mu_1f(2g-2-d_1)\geq \mu_1f(2g-2r_1)=\frac{(2g-2r_1+1)(2g+2r_1-4)+2g-2r_1}{2g-2r_1+1}\mu_1 \geq \\
&\geq \frac{(2g+2r_1-4)(2g-2r_1+1)+2g-2r_1}{2g-2r_1+1}\frac{\chi_f}{g-u_f},
\end{split}
\end{equation*}
as wanted.

Let us now suppose that the equality holds in \eqref{eq: r1}.
 If $r_1=1$, we proved in \Cref{teo: konnoplus} that necessarily $u_f=g-1$, and so $f$ is a trigonal fibration by  \Cref{prop: estremale} below. 
 If $r_2\geq 2$, then the fibre of $\cE_1$ induces a linear series on $F$, of degree $d_1=2r_1-2$. 
 If $r_1<g$, again Clifford's Theorem tells us that $F$ is hyperelliptic. So we can use the results of \cite{L-Z-hyper}, in particular the bound of theorem 1.4, that is always strictly greater than the one of \eqref{eq: r1}, a contradiction.
 If $r_1=g$, then the Hodge bundle is semistable, and in particular $u_f=0$. The equality in  \eqref{eq: r1} becomes the equality in the classical slope inequality \eqref{eq: canslope}
 \begin{equation}\label{eq: c}
 K_f^2=\frac{4g-4}{g}\chi_f.
 \end{equation}
We can use  the inequality \eqref{eq: cinque} of Konno, whose coefficient is strictly greater than the one of \eqref{eq: c}, unless $g=2$.
\end{proof}
\begin{remark}
The function bounding the slope in inequality \eqref{eq: r1} is increasing with $u_f$ and $r_1$.
For $u_f$ and $r_1$ bot small these inequalities are not meaningful. For example for the smallest cases $u_f=0$ and $r_1=1$ it is worse than the slope inequality. 
Notice that in the case where $r_1=1$, we obtain again precisely the bound \eqref{eq: konnoplus} of  \Cref{teo: konnoplus}.
Note moreover that $r_1\leq g-u_f$, so the best value for $r_1$ is if it is precisely $g-u_f$. This happens under the assumption of \Cref{teo: eurekina}.
The strength of \eqref{eq: r1} lies mainly in the case where $u_f$ is very big.
\end{remark}
We can somehow  improve this inequality by using a sharpening of Clifford's theorem obtained by the author with Riva in \cite[Theorem 2.13]{R-S}. 
Let us recall the result in the form needed here: 

\begin{theorem}[Riva-Stoppino]\label{teo: rs}
Let $F\subseteq \mathbb{P}^{g-1}$ be a canonical non-hyperelliptic curve of genus $g\geq 2$. For any subspace $W \subset H^0(K_F)$  of dimension $\dim W \geq 2$ and codimension $k$, we have:
$$\frac{\deg |W|}{\dim |W|}\geq \frac{2g-2-\beta}{g-\beta-1},$$
where $\beta:=\min \{ k, \cliff(C) \}$.
\end{theorem}

\begin{theorem}\label{teo: tecnico plus}
Let $f\colon S\to B$ a relatively minimal fibred surface of genus $g\geq 2$. Let $m:=min\{c_f, g-r_1\}$.
The following inequality holds:

\begin{equation}\label{eq: r1 plus}
K^2_f\geq \left(2g-1+\frac{2g-2-m}{g-m-1}(r_1-1)-\frac{g-m-1}{(2g-1)(g-m-1)-(2g-2-m)(r_1-1)}\right)\frac{\chi_f}{g-u_f}.
\end{equation}
\end{theorem}
\begin{proof}
Let us start two simple remarks.
(1) In the case where the fibration is hyperelliptic (i.e., $c_f=0$), we have $m=0$ and  the  inequality \eqref{eq: r1 plus} becomes precisely the inequality \eqref{eq: r1} of the theorem above. 
(2) If $r_1=1$, then the inequality \eqref{eq: r1 plus}  becomes the inequality  \eqref{eq: konnoplus}. 
 
Therefore, we have reduced ourselves to the case where $F$ is non-hyperelliptic and $r_1\geq 2$. 
In this case we can apply  \Cref{teo: rs}.  
Observe that 
\[\alpha\leq a\leq 2g-2-\frac{2g-2-m}{g-m-1}(r_1-1).\]
Now the statement follows by applying \Cref{lem: konnoplus}.
 \end{proof}

In the case where $r_1\geq 2$ we can give sharper inequalities as follows:
\begin{theorem}\label{teo: r>1}
Let $f\colon S\to B$ a relatively minimal fibred surface of genus $g\geq 2$. 
Let us suppose that $r_1=\rk\cE_1\geq 2$. Consider  the special divisor $R_1:=(K_f-Z_1)_{|F}$ on the general fibre $F$. 

Then one of the following possibilities holds:
\begin{itemize}
\item[(i)] If $h^1(R_1)\geq 2$, then 
\begin{equation}\label{eq: mah}
\begin{split}
K_f^2&\geq \left(\frac{2g+2r_1+c_f-4}{g-u_f}+\frac{2g-2r_1-c_f}{(2g-2r_1-c_f+1)(g-u_f)}\right)\chi_f\geq\\
&\geq \left(\frac{2g+c_f}{g-u_f}+\frac{2g-4-c_f}{(2g-3-c_f)(g-u_f)}\right)\chi_f.\\
\end{split}
\end{equation}
\item[(ii)] If $h^1(R_1)=1$, then necessarily $R_1\not =0$
 and the following inequality holds:
\begin{equation}\label{eq: meh}
K_f^2\geq \left(\frac{3g+r_1-4}{g-u_f}+\frac{g-r_1}{(g-r_1+1)(g-u_f)}\right)\chi_f\geq \left(\frac{3g-1}{g-u_f}+\frac{g-2}{(g-1)(g-u_f)}\right)\chi_f.
\end{equation}
\end{itemize}
\end{theorem}
\begin{proof}
Observe that $R_1$ is a special divisor of degree $d_1$ and ${M_1}_{|F}\subseteq H^0(F,R_1)$, and so $h^0(R_1)\geq r_1\geq 2$ by assumption. 

(i) If $h^1(R_1)\geq 2$, then $R_1$ contributes to the Clifford index $c_f$ of $F$ and so $d_1\geq 2(h^0(R_1)-1)+c_f\geq 2r_1-2+c_f$. So we have that $a_1\leq 2g-2r_1-c_f$, and we obtain equation \eqref{eq: mah} by using \eqref{eq: 2}.

(ii) If $h^1(R_1)=1$, in theory $R_1$ could be $0$, but in this case we would have necessarily $d_1=0$ and so $r_1=1$, contrary to the assumption. 
Therefore, necessarily $R_1\not=0$, and by Riemann-Roch \[d_1-g+2=h^0(D_1)\geq r_1,\] hence $a_1\leq g-r_1$, and we just use this bound again in inequality \eqref{eq: 2}.
\end{proof}

\subsection{Slope inequalities for Harder-Narasimhan filtration of length $\leq 2$}

Now we prove some new slope inequalities for the case where the Harder-Narasimhan filtration is short.
For the case $\ell=1$, i.e., $f_*\omega_f$ is semistable, recall the following beautiful result:

\begin{theorem}[Konno \cite{konno3}]
Let $f\colon S\to B$ a relatively minimal non locally trivial fibred surface of genus $g\geq 2$.
If the Hodge bundle is semistable, then the following inequality holds:
\begin{equation}\label{eq: cinque}
K_f^2\geq \frac{5g-6}{g}\chi_f.
\end{equation}
\end{theorem}
\begin{proof}
We give the proof because it is very short: being $f_*\omega_f$ semistable, so is its symmetric product $\sym^2 f_*\omega_f$. 
Consider the relative multiplication morphism $\gamma \colon \sym^2 f_*\omega_f \to f_*\omega_f^{\otimes 2}$. By Noether's theorem, $im\gamma$ has the same rank of $f_*\omega_f^{\otimes 2}$, so
 \[\mu(\sym^2 f_*\omega_f)\leq \mu(im\gamma)\leq \mu(f_*\omega_f^{\otimes 2}).\]
 But $\deg f_*\omega_f^{\otimes 2}= K_f^2+\chi_f$ (by Grothendieck-Riemann-Roch), and $\mu(\sym^2 f_*\omega_f)=2\mu(f_*\omega_f)=2\chi_f/g$.
 So, the inequality above becomes 
 \[2\frac{\chi_f}{g}\leq \frac{K_f^2+\chi_f}{3g-3},\]
 which gives the statement.
\end{proof}
\begin{remark}
The result above gives a high bound for the case in which the Harder-Narasimhan filtration is of length $1$. 
It is not known  whether  this bound is sharp or not.
\end{remark}

We now give a new bound for the slope under the assumption that  Harder-Narasimhan filtration of $\cE=f_*\omega_f$ has length $ 2$ and  that $\mu_l=0$ (i.e., $u_f>0$). 
\begin{theorem}\label{teo: eurekina}
Let $f\colon S\to B$ be  a relatively minimal  locally non-trivial fibred surface of genus $g\geq 2$, such that $\cA$ is semistable and $\mu_l=0$.
Then the following inequality holds:
\begin{equation}\label{eq: eureka}
K_f^2\geq \left( \frac{4(g-1)-2u_f}{(g-u_f)}+ \frac{2u_f}{(2u_f+1)(g-u_f)}\right)\chi_f.
\end{equation}
Moreover, assume that the  equality holds in \eqref{eq: eureka}.
Then necessarily
$u_f=g-1$ and  $F$ is trigonal.

Assume that $u_f\leq c_f$. 
Then the following inequality holds:
\begin{equation}\label{eq: eureka-bis}
K_f^2\geq \left( \frac{4(g-1)-u_f}{(g-u_f)}+ \frac{u_f}{(u_f+1)(g-u_f)}\right)\chi_f.
\end{equation}
Assume that $u_f\geq c_f$. 
Then the following inequality holds:
\begin{equation}\label{eq: eureka-biss}
K_f^2\geq \left( 2g-1+\frac{2g-2-c_f}{g-c_f-1}(g-u_f-1) -\frac{g-c_f-1}{2gu_f-gc_f-2u_f-c_fu_f-1}\right)\frac{\chi_f}{(g-u_f)}.
\end{equation}
\end{theorem}
\begin{proof}
From the assumptions we have that the  Harder-Narasimhan filtration of the Hodge bundle $\cE$ has length 2:
\[0=\cE_0\subset \cE_1=\cA\subset \cE_2=\cE.\] 
Moreover, $\mu_1=\chi_f/(g-u_f)$ and $\mu_2=0$.
Observe that under our assumptions $r_1=g-u_f$.

We use  \Cref{teo: tecnico}: the inequality \eqref{eq: r1} becomes 
\[K^2_f\geq \left(\frac{4g-4-2u_f}{g-u_f}+\frac{2u_f}{(2u_f+1)(g-u_f)}\right)\chi_f,\]
as wanted.

The case where the equality holds follows again from \Cref{teo: tecnico} (recalling that in this case $r_1<g$).
As for the other two inequalities, they follow directly from  \Cref{teo: tecnico plus}.
\end{proof}

\begin{remark}
As explained for example in \cite[Remark 3.3]{SB}, using Xiao's method with the assumptions of \Cref{teo: eurekina} one obtains
the worse inequality:
\[K_f^2\geq \left( \frac{4(g-1)-2u_f}{(g-u_f)}\right)\chi_f.\]
This very same inequality is reproved in the article \cite{CLZ} in the proof of Lemma 3.3.2.
\end{remark}


\section{A discussion on the extremal cases}\label{sec: extremal}

We now list some results, due to Serrano, Pirola, Zucconi and Beorchia, about locally non-trivial fibred surfaces with maximal unitary rank $u_f=g-1$

First of all, in \cite[Theorem 3.6 and Remark 3.7]{pietroBUMI} Pirola proved that in this case up to an \'etale base change $g-1=u_f=q_f$ (see  \Cref{rem: confronto} to see that this is not obvious).

\begin{proposition}\label{prop: estremale2}
Let $f\colon S\to B$ a relatively minimal non-trivial fibred surface of genus $g\geq 4$ such that $u_f= g-1$. 
Then the general fibre $F$ is trigonal;
\end{proposition}
\begin{proof}
 By Pirola's result we can assume that the relative irregularity is maximal.
Cai proved in \cite{Cai} that for locally non-trivial hyperelliptic  fibrations Xiao's bound \eqref{eq: xiao bound} $q_f\leq (g+1)/2$ holds, so for $g\geq 3$ no hyperelliptic fibration can have $q_f=g-1$. 
In the case where $g=2$ of course $F$ is hyperelliptic. The case $g=3$ has been  investigated in \cite{moller}. 
On the other hand, by  Barja-Naranjo-Gonzalez \cite{BGN}, we have that $q_f\leq g-c_f$. 
Moreover, Favale-Naranjo-Pirola \cite{FNP} proved that if the general fibre is a degree $d$ plane curve (thus of gonality $d-1$ and Clifford index $d-3$), then $q_f\leq g-c_f-1$. 
Therefore, under the assumption of the proposition, $\cliff(F)\leq 1$ and $F$ is not a plane quintic, nor it is hyperelliptic (i.e., with $\cliff(F)>0$). Hence we conclude that $F$ necessarily is trigonal. 
\end{proof}
\begin{theorem}\label{teo: estremale}
Let $f\colon S\to B$ be a relatively minimal non-trivial fibred surface of genus $g\geq 4$ such that $u_f= g-1$. 
Then the following hold:
\begin{enumerate}
\item[(1)]  the fibration is non-isotrivial;
\item[(2)] the smooth fibres $F_t$ are a covering of a non-isotrivial family of elliptic curves $E_t$.
\item[(3)] the fibration is not Kodaira , i.e. $K_f^2<12\chi_f$.
\end{enumerate}
\end{theorem}

\begin{proof}
 By Pirola's result we can assume that after an \'etale base change, $q_f=u_f=g-1$. \'Etale base changes do not modify the slope.
 
(1) This is due to Serrano: indeed, in \cite{serrano} he proved that Xiao's bound \eqref{eq: xiao bound} $q_f\leq (g+1)/2$ holds for isotrivial fibrations. So we can exclude the isotrivial fibrations unless $g\leq 3$. In \cite{moller} isotrivial fibrations of genus $2$ and $3$ with irregularity respectively $1$ and $2$ are constructed. 

(2) This proof is due to Pietro Pirola in a private communication. 
Let $B^0$  be the smooth locus of $f$ (i.e., the locus consisting of all $t\in B$ such that the fibre $F_t$ is smooth). 
For any $t\in B^0$, let $j_t\colon F_t\hookrightarrow S$ be the inclusion.
We have the following diagram 
\[
\xymatrix{
F_t\ar[d]\ar@{^{(}->}^{j_t}[r]&S \ar^{\alb_S}[d]\\
J(F_t)\ar[r]^{{j_t}_*} &\Alb(S)}
\]
Now, the morphism ${j_t}_*$ goes in a family of translates of a fixed abelian sub-variety $A\subset \Alb(S)$ of dimension $q_f=u_f=g-1$.
Let us consider the induced surjective homomorphism $g_t\colon J(F_t)\to A$, and let be $(\ker(g_t)_0)$ be the abelian subvarieties of $J(F_t)$ that are  the connected component of $\ker(g_t)$ containing the origin, for any $t\in B^0$. Note that $(\ker(g_t)_0)$ is an elliptic curve for any $t\in B^0$.
Passing to the dual varieties and homomorphisms we have surjective homomorphisms
\[ \pic^0(F_t)\to (\ker(g_t)_0)^\vee=:E_t\] for any $t\in B^0$. 
We thus have the morphisms
\[F_t\hookrightarrow  \pic^0(F_t) \longrightarrow E_t.\]
The composition $F_t\to E_t$ necessarily   is  surjective because the image of $F_t$ generates $\pic^0(F_t)$ as an abelian variety.
Therefore, this is the desired morphism of the general fibre to a family of elliptic curves. 

Now, the family  $E_t$ has varying moduli: if this were not the case then all $\pic^0(F_t)$ would be a fixed abelian variety.  
This would imply that  the entire family over $B_0$ would be trivial up to a base change, a contradiction by (1).

(3) Also this  proof is due to Pietro Pirola in a private communication. By performing  a base change $B'\to B$ if necessary, we have a global non-isotrivial family of genus $1$ curves $W\stackrel{\alpha}\to B'$ such that, calling $f'\colon S'\to B'$ the fibration obtained by base change, there exists a covering 
$\beta \colon S'\to W$ such that $f'=\alpha\circ\beta$. 
\[
\xymatrix{
S'\ar_{f'}[dd]\ar^{\beta}[dr]\ar[rr]& &S \ar^{f}[dd]\\
& W\ar^\alpha[dl] &\\
B'\ar[rr]& & B
}
\]
Now the fibration of genus one contains at least one fibre with non compact Jacobian. Thus $f$ contains at least a fibre with the same property, and this implies that it is not isotrivial.
\end{proof}

\section{An upper bound on the unitary rank}\label{sec: max irr}

 We now extend Konno's bound \cite[Lemma 2.7]{konno} on $q_f$ to $u_f$, and improve it to a strict inequality.
 \begin{theorem}\label{thm: bound}
 Let $f\colon S\to B$ be a relatively minimal non-isotrivial fibration of genus $g\geq 2$. Then
\begin{equation}\label{eq: xiaoplus}
u_f < g \frac{5g-2}{6g-3}.
\end{equation}
\end{theorem}
\begin{proof}
The non-sharp inequality is immediate by combining the inequalities \eqref{eq: konnoplus} and \eqref{eq: 12}.
Suppose now that the equality holds in \eqref{eq: xiaoplus}. 
Then the equality in \eqref{eq: konnoplus} holds, so $f$ is necessarily of maximal unitary rank $u_f=g-1$, and the equality in \eqref{eq: 12} holds; so  $f$ is necessarily Kodaira. 
We have seen in \Cref{teo: estremale} that this case cannot happen.
\end{proof}

\begin{theorem}
Let $f\colon S\to B$ be a relatively minimal non-isotrivial fibration of genus $g\geq 2$. 
Suppose that $u_f$ is maximal, i.e., $u_f=g-1$.
Then $g\leq 6$.
\end{theorem}
\begin{proof}
The inequality \eqref{eq: konnoplus} in the first case becomes 
\[K_f^2\geq \frac{4g(g-1)}{2g-1}\chi_f.\]
If we combine it with the bound \eqref{eq: 12}, we obtain $4g(g-1)\leq 3(2g-1)$, so $g^2-7g+3\leq 0$, which gives that $g\leq 6$.
\end{proof}
We can also give a bound in the case $u_f=g-2$. 
By using the inequality \eqref{eq: xiaoplus} in this case we obtain 
the equation 
\[g^2-13g+7\leq 0,\]
which implies
$g\leq 12$. 
We can rule out the case where $g=12$ as follows.
\begin{proposition}\label{prop: estremale}
Let $f\colon S\to B$ a relatively minimal non locally trivial fibred surface of genus $g\geq 2$. 

Suppose that $u_f= g-2$. 
Then the genus of the fibration is $g\leq 11$.  
\end{proposition}
\begin{proof}
The inequality \eqref{eq: konnoplus} becomes 
\[K_f^2\geq \frac{2g(g-1)}{2g-1}\chi_f,\]
and so, combining again with \eqref{eq: 12} we get $g\leq 12$. 
But for $g=12$ we obtain that the equality is reached in \eqref{eq: konnoplus} of \Cref{teo: konnoplus}, so $u_f$ should be $g-1$, a contradiction.
\end{proof}
I the paper  \cite{LS2} we give many conditions on the case where $g=6$ and $u_f=g-1=5$.


\section{An application to the Coleman-Oort conjecture}\label{sec: CO}
We will use the setting and the notations of the papers \cite{LZ4} and \cite {CLZ}. Let us consider $\cM_g=\cM_{g,\ell}$ the moduli space of smooth projective curves of genus $g\geq 2 $ with a full level $\ell-$structure. and $\cA_g=\cA_{g,\ell}$ the moduli space of $g$-dimensional principally polarized abelian varieties with a full level $\ell-$structure. The level $\ell\geq3$ assures that these are moduli schemes representing the corresponding moduli functor. 
Let 
\[
j^\circ\colon \cM_g\longrightarrow \cA_g
\]
be the Torelli morphism. The Torelli locus $\mathcal T_g$ is the closure in $\mathcal A_g$ of the image of $j^\circ$, which we call $\mathcal T_g^\circ$.

A closed positive-dimensional subvariety $Z\subseteq \cA_g$ is said to be {\em contained generically in the Torelli locus} if $Z\subset \cT_g$ and $Z\cap \cT_g^\circ\not=\emptyset$.

The moduli space $\cA_g$ is isomorphic to a Shimura variety. In $\cA_g$ there are Shimura subvarieties. 
The famous Coleman-Oort conjecture predicts that for large $g$ there does not exist any positive-dimensional Shimura subvariety generically contained in the Torelli locus. 

For $g\leq 4$ such subvarieties do exist: see \cite{FGS}.
Many authors have proved results in the direction of this conjecture, let us recall at least 
\cite{FGP}, \cite{CFG}, \cite{GP}, \cite{LZ4}.

Recall that, given $C\subset \mathcal A_g$ a smooth closed curve, the canonical Higgs bundle $\mathcal E_C$ on $C$ is the Hodge bundle given by the universal family of abelian varieties restricted over $C$. The bundle $\mathcal E_C$ decomposes as $\mathcal E_C=\mathcal F_C\oplus \mathcal U_C$ where $\mathcal U_C$ is the maximal unitary Higgs subbundle corresponding to the maximal sub-representation on which $\pi_1(C)$ acts through a compact unitary group. The bundles $\mathcal E_C$, $\mathcal F_C$ and $\mathcal U_C$ have an Hodge decomposition into a $(-1,0)$ and $(0,-1)$ part, which we call $\cU_C^{-1,0 }$ and $\cU_C^{0,-1}$.

Here, using the results of the previous sections,  we can improve  a  result of Chen,  Lu and Zuo in \cite[Thm 1.1.2]{CLZ}.
We give a condition on the rank of the $(-1,0)$ part of  the maximal unitary Higgs subbundle $\mathcal U_C$ of a curve $C\subset \mathcal A_g$ that implies that it is not generically contained in the Torelli locus.
\begin{theorem}\label{teo: CO}
Let $C \subset \cA_g$ be a curve with Higgs bundle decomposition $\cE_C = \cA_C \oplus \cU_C$, where $\cU_C$ is the maximal unitary Higgs subbundle. 
\begin{itemize}
\item[(i)] If
\begin{equation}\rk \,\cU_C^{-1,0 }\geq g \frac{5g-2}{6g-3},
\end{equation} then $C$ is not contained generically in the Torelli locus.
\item[(ii)]
If $C$ is a Shimura curve such that 
\begin{equation}
\rk\, \cU_C^{-1,0} \geq \frac{4g-1+\sqrt{16g^2-36g+21}}{10},
\end{equation}
 then $C$ is not contained generically in the Torelli locus.
\end{itemize}
\end{theorem}
\begin{proof}
The proof is completely analogous to the one of \cite[Thm 1.12]{CLZ}. The key point is that if $C$ is contained generically in the Torelli locus, then we can consider 
 the following diagram
$$
\xymatrix{
B\ar[r]^{j_{B}}\ar[d]& C\ar[d]\\
\cM_g \ar[r]^{j^\circ}&\cA_g\\
}
$$
where $B$ is the normalization of the pullback $(j^\circ)^{-1}C$ in $\cM_g$. The morphism $B\to \cM_g$ gives a fibration $f\colon S\to B$  induced by the fine moduli structure. We can perform the semistable reduction to a compactification of the family  obtaining a  non-isotrivial fibred surface $\overline f\colon \overline S\to \overline B$, and a morphism $\overline B\to \overline {\cM}_g$ extending $B\to \cM_g$.

Now, the Hodge bundle $\overline f_*\omega_{\overline f}$  is the $(-1,0)-$part of the Higgs bundle on $\overline B$ associated to the relative semistable Jacobian family (see \cite[Sec. 3]{LZ-arxive}). The pullback via $j_{\overline B}$ of the $(-1,0)$ part of the decomposition 
$\mathcal E_{\overline B}=\mathcal F_{\overline B}\oplus \mathcal U_{\overline B}$ is precisely the 
second Fujita decomposition \eqref{second fujita} of ${\overline f}_*\omega_{\overline f}$, so 
$$\rk \, {\mathcal U_{\overline B}}^{-1,0}=u_{\overline f}.$$
The statement $(i)$ directly follows from \Cref{thm: bound}. 

In order to prove  $(ii)$, we need  a result due to Viehweg and Zuo (\cite{VZ}, see also \cite[Thm. 1.2.1]{CLZ}): if $C$ is Shimura, then $\mathcal F_{\overline B}^{-1,0}$ is semistable. 
We thus are in the assumptions of  \Cref{teo: eurekina}.
The inequality \eqref{eq: eureka} combined with \eqref{eq: 12} gives: 
\begin{align*}
& \frac{(4(g-1)-2u_{\overline f})(2u_{\overline f}+1)+2u_{\overline f}}{(2u_{\overline f}+1)(g-u_{\overline f})}\leq 12\iff\\
&(2u_{\overline f}+1)(g-1)-u_{\overline f}^2\leq3(g-u_{\overline f})(2u_{\overline f}+1)\iff\\
&5u_{\overline f}^2-(4g-1)u_{\overline f}-(2g+1)\leq 0.
\end{align*}
Now, the discriminant of this quadratic form is $\Delta=16g^2-36g+21>0$, and so we need to have 
\begin{equation}\label{ye}
u_{\overline f}\leq \frac{4g-1+\sqrt{16g^2-36g+21}}{10}.
\end{equation}
If the equality holds, then in particular $f$ is Kodaira, and the equality holds in \eqref{eq: eureka}, so by \Cref{teo: eurekina} we have that:
$u_{\overline f}=g-1$. In this case we already proved in  \Cref{teo: estremale} that the fibration is not Kodaira.
\end{proof}
\begin{remark}\label{rem: confronto LZ}
The bounds found by Chen, Lu and Zuo in \cite[Thm 1.12]{CLZ} are respectively:
\begin{itemize}
\item[(i)] if $\rk \,\cU_C^{-1,0 }>  \frac{5g+1}{6},$ then $C$ is not generically contained in the Torelli locus;
\item[(ii)] if $\rk \,\cU_C^{-1,0} \geq \frac{4g+2}{5}$ and $C$ is Shimura, then it is not generically contained in the Torelli locus.
\end{itemize} 
So, as $g \frac{5g-2}{6g-3}<  \frac{5g+1}{6}$ and $\frac{4g-1+\sqrt{16g^2-36g+21}}{10}<\frac{4g-3}{5}< \frac{4g+2}{5}$, we see that  we are giving a small improvement to their results. 
\end{remark}


\noindent {Lidia Stoppino,\\Dipartimento di Matematica, \\Universit\`a di Pavia,\\ Via Ferrata 5, 27100 Pavia, Italy.\\
E-mail: \textsl {lidia.stoppino@unipv.it}

\end{document}